\documentclass[11pt,reqno]{amsart}
\usepackage{amsmath,amssymb}
\usepackage{hyperref}
\usepackage{url}
\usepackage{tikz,enumerate}
\usepackage{diagbox}
\usepackage{appendix}
\usepackage{epic}
\allowdisplaybreaks[4]

\theoremstyle{definition}

\theoremstyle{remark}

\numberwithin{equation}{section}
\numberwithin{theorem}{section}
\numberwithin{defn}{section}

\begin{document}
\title[New Proofs of Ramanujan's Identities on False Theta Functions]
 {New Proofs of Ramanujan's Identities on \\ False Theta Functions}

\author{Liuquan Wang}
\address{School of Mathematics and Statistics, Wuhan University, Wuhan 430072, Hubei, People's Republic of China}
\email{mathlqwang@163.com;wanglq@whu.edu.cn}



\dedicatory{}

\maketitle
\begin{abstract}
We provide new proofs to five of Ramanujan's intriguing identities on false theta functions without using the Rogers-Fine identity and Bailey transforms.
\end{abstract}

\section{Introduction}

In 1917, Rogers \cite{rogers} introduced false theta functions, which are series that are instances of classical theta series except for an alteration of the signs of some of the series' terms. False theta series behave quite different from theta series and their properties are less well understood.  In his notebooks \cite{notebook} and lost notebook \cite{lostnotebook}, Ramanujan studied several examples of false theta function identities.
For example, on page 13 of Ramanujan's Lost Notebook \cite{lostnotebook} (c.f. \cite[Section 9.3, pp. 227--232]{lost-notebook1}), the following five identities were recorded:
\begin{align}
\sum_{n=0}^{\infty} \frac{(-1)^n (q;q^2)_n \, q^{n(n+1)}}{(-q;q)_{2n+1}} & =\sum_{n=0}^{\infty} (-1)^n q^{n(n+1)/2}, \label{false-1}\\
\sum_{n=0}^{\infty} \frac{(q;q^2)_n^2 \, q^{n} }{(-q;q)_{2n+1}} & =\sum_{n=0}^{\infty} (-1)^n q^{n(n+1)}, \label{false-2}\\
\sum_{n=0}^{\infty} \frac{(q;q^2)_n \, q^n }{(-q;q)_{2n+1}} & =\sum_{n=0}^{\infty} (-1)^n q^{3n(n+1)/2},\label{3} \\
\sum_{n=0}^{\infty} \frac{(q;-q)_{2n} \, q^n}{(-q;q)_{2n+1}} & =\sum_{n=0}^{\infty} (-1)^n q^{2n(n+1)}, \label{false-4}\\
\sum_{n=0}^{\infty} \frac{(q;-q)_n (-q^2;q^2)_n \, q^{n}}{(-q;q)_{2n+1}} & =\sum_{n=0}^{\infty} (-1)^n q^{3n(n+1)}, \label{false-5}
\end{align}
where
\[(a;q)_{n}=(1-a)(1-aq)\cdots (1-aq^{n-1}).\]
Note that the right hand sides of these identities are instances of series
\begin{equation*}
\sum_{n=0}^{\infty} (-1)^n q^{n(n+1)/2}.
\end{equation*}

Ramanujan's false theta function identities \eqref{false-1}--\eqref{false-5} were first proved by Andrews \cite[Section 6]{Andrews-Adv} by employing the Rogers-Fine identity as well as some other identities. Recently, Andrews and Warnaar \cite{Andrews-Warnaar} provided new proofs of \eqref{false-1}--\eqref{false-4} by using symmetric bilateral Bailey transforms. However, their methods apparently do not yield a proof of \eqref{false-5}. Using other bilateral Bailey transformation formulaes, Chu and Zhang \cite{Chu-Zhang} proved many false theta function identities including \eqref{false-1}--\eqref{false-5}.

Inspired by their work, the goal of this paper is to provide new proofs of the identities \eqref{false-1}--\eqref{false-5} without using the Rogers-Fine identity or Bailey transform. Instead, the key tool in our proof is a beautiful formula given by Liu \cite{Liu}. Our proofs of \eqref{false-1}--\eqref{false-5} share the same pattern and are motivated by a recent work of the present author and Yee \cite{Wang-Yee}. Our method is constructive, which means that we can get the right side from the left side of an identity directly without knowing the right side in advance. We believe that the method used here can be applied to discover and prove many other identities.

\section{Preliminaries}\label{sec2}
In this section, we collect some identities on basic hypergeometric series in the literature for later use.

We use the following standard $q$-series notation:
\begin{align*}
(a;q)_{\infty}&:= \lim_{n\to \infty} (a;q)_n, \quad |q|<1, \\
(a_1,a_2,\dots, a_r;q)_{\infty}&:=(a_1;q)_{\infty}(a_2;q)_{\infty}\cdots (a_r;q)_{\infty}
\end{align*}
and
\begin{equation*}\label{bhs}
{}_{r+1}\phi_{r}\bigg(\genfrac{}{}{0pt}{}{a_0,  a_1, \dots,  a_r}
  {b_1,   \dots,  b_{r}}; q,
z \bigg) :=\sum_{n=0}^{\infty} \frac{(a_0;q)_n (a_1;q)_n \cdots (a_r;q)_n}{(q;q)_n (b_1;q)_n \cdots (b_{r};q)_n} z^n.
\end{equation*}

Firstly, we need some ${}_{2}\phi_{1}$ formulas. One of Heine's famous transformation formula for ${_2}\phi_{1}$ series (see \cite[p.13]{gasper}) is
\begin{align}\label{Heine}
{_2}\phi_{1}\bigg(\genfrac{}{}{0pt}{} {a,  b}{c};q,z \bigg)=\frac{(b,az;q)_{\infty}}{(c,z;q)_{\infty}}{_2}\phi_{1}\bigg(\genfrac{}{}{0pt}{} {c/b,  z} {az}; q,b\bigg).
\end{align}
From \cite[p.\ 15]{gasper}, we also find that
\begin{equation}
{_2} \phi_1 \bigg(\genfrac{}{}{0pt}{} {a,  b} {c}; q, z\bigg) =\frac{(az;q)_{\infty}}{(z;q)_{\infty}} \sum_{n=0}^{\infty} \frac{(a, c/b;q)_{n}}{(q,c,az;q)_n} (-bz)^n q^{\binom{n}{2}}. \label{gasperp15}
\end{equation}

Another useful formula can be found in \cite[p.71]{gasper}:
\begin{align}\label{gasper71-new}
{_2}\phi_{1}\bigg(\genfrac{}{}{0pt}{} {q^{-n}, d/b} {d};q,bq/e\bigg)=(-1)^{n}q^{-\binom{n}{2}}(e;q)_{n}e^{-n}{_3}\phi_{2}\bigg(\genfrac{}{}{0pt}{} {q^{-n}, b, 0}{d, e};q,q\bigg).
\end{align}

The following identity is due to Bailey and Daum \cite[Eq.\ (1.8.1)]{gasper}:
\begin{align}\label{BD}
{_2}\phi_{1}\bigg(\genfrac{}{}{0pt}{} {a,  b} {aq/b}; q,-\frac{q}{b} \bigg)=\frac{(-q;q)_{\infty}(qa,aq^2/b^2;q^2)_{\infty}}{(aq/b,-q/b;q)_{\infty}}.
\end{align}

Next, we also need the  $q$-Pfaff--Saalsch\"utz summation formula \cite[p.40, Eq. (2.2.1)]{gasper}:
\begin{align}
{_3}\phi_2 \bigg(\genfrac{}{}{0pt}{} {q^{-n},  aq^n,  aq/bc} {aq/b,  aq/c}; q,  q \bigg) &= \frac{(b,c;q)_n}{(aq/b, aq/c ;q)_n}  \left(\frac{aq}{bc}\right)^n.  \qquad \label{gasper221}
\end{align}

Our proofs of \eqref{false-1}--\eqref{false-5} are all based on the following beautiful transformation formula due to Liu \cite{Liu}: for $|\alpha a b/q|<1$,
\begin{align}
&\frac{(\alpha q, \alpha a b/q;q)_{\infty}}{(\alpha a, \alpha b;q)_{\infty}} {_3}\phi_2 \bigg(\genfrac{}{}{0pt}{} {q/a, q/b,  \beta} {c, d}; q , \frac{\alpha a b}{q} \bigg) \notag\\
&=\sum_{n=0}^{\infty} \frac{(1-\alpha q^{2n}) (\alpha, q/a, q/b;q)_n (-\alpha ab/q)^n q^{\binom{n}{2}}}{(1-\alpha) (q, \alpha a, \alpha b;q)_n} {_3}\phi_2 \bigg(\genfrac{}{}{0pt}{} {q^{-n}, \alpha q^n , \beta}  {c, d} ; q , q  \bigg). \label{liu2089}
\end{align}

\section{Proofs of Ramanujan's Identities}
Now we give new proofs to \eqref{false-1}--\eqref{false-5}.
\begin{proof}[Proof of \eqref{false-1}]
We observe that
\begin{align} \label{1-start}
&\sum_{n=0}^{\infty}\frac{(-1)^{n}(q;q^2)_{n}q^{n(n+1)}}{(-q^2;q^2)_{n}(-q^3;q^2)_{n}} \nonumber \\
=& \lim_{a \to 0}{_3}\phi_{2}\bigg(\genfrac{}{}{0pt}{} {q, q^2,  q^2/a} {-q^2, -q^3}; q^2, a \bigg).
\end{align}
Replacing $q$ by $q^2$ and setting $(b,c,d,\alpha,\beta)=(1,-q^2,-q^3,q^2,q)$ in \eqref{liu2089}, we get
\begin{align}\label{1-middle}
&\frac{(q^4,a;q^2)_{\infty}}{(q^2a,q^2;q^2)_{\infty}}{_3}\phi_{2}\bigg(\genfrac{}{}{0pt}{} {q, q^2,  q^2/a} {-q^2, -q^3}; q^2, a \bigg) \nonumber\\
=&\sum_{n=0}^{\infty}\frac{(1-q^{4n+2})(q^2,q^2/a,q^2;q^2)_{n}(-a)^{n}q^{n(n-1)}}{(1-q^2)(q^2,q^2a,q^2;q^2)_{n}} {_3}\phi_{2}\bigg(\genfrac{}{}{0pt}{} {q^{-2n}, q^{2n+2},q}  {-q^2, -q^3};q^2,q^2 \bigg).
\end{align}
Substituting this into \eqref{1-start} and taking the $a\to 0$ limit, we obtain
\begin{align}\label{1-final}
&\sum_{n=0}^{\infty}\frac{(-1)^{n}(q;q^2)_{n}q^{n(n+1)}}{(-q^2;q^2)_{n}(-q^3;q^2)_{n}} \nonumber \\
=& \sum_{n=0}^{\infty}(1-q^{4n+2})q^{2n^2}{_3}\phi_{2}\bigg(\genfrac{}{}{0pt}{} {q^{-2n},  q^{2n+2}, q} {-q^2, -q^3};q^2,q^2 \bigg).
\end{align}
The ${}_{3}\phi_{2}$ series in the right sum equals $\frac{(1+q)q^n}{1+q^{2n+1}}$ by the $q$-Pfaff--Saalsch\"utz summation, so that we get
\begin{align*}
\sum_{n=0}^{\infty}\frac{(-1)^{n}(q;q^2)_{n}q^{n(n+1)}}{(-q;q)_{2n+1}} =&\sum_{n=0}^{\infty}(1-q^{2n+1})q^{2n^2+n} \nonumber \\
=&\sum_{n=0}^{\infty}q^{2n^2+n}-\sum_{n=0}^{\infty}q^{(2n+1)(n+1)} \nonumber \\
=&\sum_{n=0}^{\infty}q^{2n(2n+1)/2}-\sum_{n=0}^{\infty}q^{(2n+1)(2n+2)/2} \nonumber \\
=&\sum_{n=0}^{\infty}(-1)^{n}q^{n(n+1)/2}. \qedhere
\end{align*}
\end{proof}


\begin{proof}[Proof of \eqref{false-2}]
We observe that
\begin{align}
\sum_{n=0}^{\infty} \frac{(q;q^2)_n^2 \, q^{n} }{(-q;q)_{2n+1}} =\frac{1}{1+q}{}_{3}\phi_{2}\bigg(\genfrac{}{}{0pt}{} {q^2, q, q}{-q^2,-q^3};q^2,q\bigg). \label{2-start-add}
\end{align}
To evaluate the ${}_{3}\phi_{2}$ series on the right side, we use \eqref{liu2089}. Replacing $q$ by $q^2$ and setting $(a,b,c,d,\alpha,\beta)=(aq,1,-q^2,-q^3,q^2/a,aq)$ in \eqref{liu2089}, we obtain
\begin{align}\label{Liu-1-coro}
&\frac{(a^{-1}q^4,q;q^2)_{\infty}}{(q^3,a^{-1}q^2;q^2)_{\infty}}{}_{3}\phi_{2}\bigg(\genfrac{}{}{0pt}{}  {a^{-1}q, q^2,  aq} {-q^2, -q^3};q^2,q \bigg) \\ =&\sum\limits_{n=0}^{\infty}\frac{(1-a^{-1}q^{4n+2})(a^{-1}q^2,a^{-1}q,q^2;q^2)_{n}(-1)^{n}q^{n^2}}{(1-a^{-1}q^2)(q^2,q^3,a^{-1}q^2;q^2)_{n}} \nonumber \\
&\cdot {}_{3}\phi_{2}\bigg(\genfrac{}{}{0pt}{} {q^{-2n}, a^{-1}q^{2+2n}, aq}{-q^2, -q^3};q^2,q^2 \bigg). \nonumber
\end{align}
Replacing $q$ by $q^2$ and taking $(a,b,c)=(a^{-1}q^2,-a^{-1}q^2,-a^{-1}q)$ in \eqref{gasper221}, we deduce that
\begin{equation}\label{value}
{}_{3}\phi_{2}\bigg(\genfrac{}{}{0pt}{} {q^{-2n}, a^{-1}q^{2+2n}, aq} {-q^2,-q^3};q^2,q^2 \bigg)=\frac{(-a^{-1}q,-a^{-1}q^2;q^2)_{n}}{(-q^2,-q^3;q^2)_{n}}(aq)^{n}.
\end{equation}
Substituting this identity into (\ref{Liu-1-coro}) and setting $a=1$, after rearranging, we get the desired identity we want.
\end{proof}

\begin{proof}[Proof of \eqref{3}]
We observe that
\begin{align}\label{3-start}
\sum_{n=0}^{\infty}\frac{(q;q^2)_{n}q^n}{(-q;q)_{2n+1}}=\frac{1}{1+q}{_3}\phi_{2}\bigg(\genfrac{}{}{0pt}{} {q^2,  q, 0}{-q^2, -q^3};q^2,q\bigg).
\end{align}
The ${}_{3}\phi_{2}$ series on the right side can be simplified in several steps. First, replacing $q$ by $q^2$ and setting $(a,b,c,d,\alpha, \beta)=(1,q,-q^2,-q^3,q^2,0)$ in \eqref{liu2089}, we deduce that
\begin{align}\label{3-continue}
&\frac{(q^4,q;q^2)_{\infty}}{(q^2,q^3;q^2)_{\infty}}{_3}\phi_{2}\bigg(\genfrac{}{}{0pt}{} {q^2,  q, 0}{-q^2, -q^3};q^2,q\bigg)\nonumber \\
&=\sum_{n=0}^{\infty}\frac{(1-q^{4n+2})(q^2,q^2,q;q^2)_{n}(-q)^{n}q^{n(n-1)}}{(1-q^2)(q^2,q^2,q^3;q^2)_{n}} \cdot {_3}\phi_{2}\bigg(\genfrac{}{}{0pt}{} {q^{-2n}, q^{2n+2}, 0}{-q^2, -q^3};q^2,q^2\bigg) \nonumber \\
&=\sum_{n=0}^{\infty}(-1)^{n}q^{n^2}\cdot \frac{1+q^{2n+1}}{1+q} \cdot {_3}\phi_{2}\bigg(\genfrac{}{}{0pt}{} {q^{-2n}, q^{2n+2}, 0}{-q^2, -q^3}; q^2,q^2\bigg).
\end{align}
Replacing $q$ by $q^2$ and setting $(b,d,e)=(q^{2n+2},-q^2,-q^3)$ in \eqref{gasper71-new}, we get
\begin{align}\label{3-gasper}
&{_2}\phi_{1}\bigg(\genfrac{}{}{0pt}{} {q^{-2n},  -q^{-2n}}{-q^2};q^2, -q^{2n+1}\bigg) \nonumber \\
&=q^{-n(n+2)}(-q^3;q^2)_{n}{_3}\phi_{2}\bigg(\genfrac{}{}{0pt}{} {q^{-2n},  q^{2n+2}, 0} {-q^2, -q^3};q^2,q^2\bigg).
\end{align}
Replacing $q$ by $q^2$ and setting $(a,b,c,z)=(q^{-2n},-q^{-2n},-q^2,-q^{2n+1})$ in \eqref{Heine}, we obtain
\begin{align}\label{3-Heine}
&{_2}\phi_{1}\bigg(\genfrac{}{}{0pt}{} {q^{-2n}, -q^{-2n}}{-q^2};q^2, -q^{2n+1}\bigg) \nonumber \\
&=\frac{(-q^{-2n},-q;q^2)_{\infty}}{(-q^2,-q^{2n+1};q^2)_{\infty}}{_2}\phi_{1}\bigg(\genfrac{}{}{0pt}{} {q^{2n+2},  -q^{2n+1}} {-q}; q^2, -q^{-2n}\bigg) \nonumber \\
&=q^{-n(n+1)}(-q;q^2)_{n}(-q^2;q^2)_{n}{_2}\phi_{1}\bigg(\genfrac{}{}{0pt}{}{-q^{2n+1}, q^{2n+2}}{-q};q^2, -q^{-2n}\bigg).
\end{align}
Replacing $q$ by $q^2$ and setting $(a,b)=(-q^{2n+1},q^{2n+2})$ in \eqref{BD}, we deduce that
\begin{align}\label{3-BD}
&{_2}\phi_{1}\bigg(\genfrac{}{}{0pt}{} {-q^{2n+1},  q^{2n+2}}{-q}; q^2, -q^{-2n}\bigg) \nonumber \\
&=\frac{(-q^2;q^2)_{\infty}(-q^{2n+3},-q^{1-2n};q^4)_{\infty}}{(-q,-q^{-2n};q^2)_{\infty}} \nonumber \\
&= \frac{q^{n(n+1)/2}}{(-q^2;q^2)_{n}}.
\end{align}
Substituting \eqref{3-BD} into \eqref{3-Heine}, we deduce that
\begin{align}\label{3-Heine-value}
{_2}\phi_{1}\bigg(\genfrac{}{}{0pt}{} {q^{-2n},  -q^{-2n}}{-q^2} ;q^2, -q^{2n+1}\bigg)=q^{-n(n+1)/2}(-q;q^2)_{n}.
\end{align}
Now substituting \eqref{3-Heine-value} into \eqref{3-gasper}, we obtain
\begin{align}\label{3-gasper-value}
{_3}\phi_{2}\bigg(\genfrac{}{}{0pt}{} {q^{-2n},  q^{2n+2},  0} {-q^2, -q^3};q^2,q^2\bigg)=\frac{(1+q)q^{(n^2+3n)/2}}{1+q^{2n+1}}.
\end{align}
Finally, substituting \eqref{3-gasper-value} into \eqref{3-continue}, we deduce that
\begin{align}
\frac{1}{1+q}{_3}\phi_{2}\bigg(\genfrac{}{}{0pt}{} {q^2,  q, 0} {-q^2, -q^3};q^2,q\bigg)=\sum_{n=0}^{\infty}(-1)^{n}q^{3n(n+1)/2}.
\end{align}
By \eqref{3-start}, we complete our proof.
\end{proof}

\begin{proof}[Proof of \eqref{false-4}]
We observe that
\begin{align}\label{4-start}
\sum_{n=0}^{\infty}\frac{(q;-q)_{2n}q^n}{(-q;q)_{2n+1}}&=\frac{1}{1+q}\sum_{n=0}^{\infty}\frac{(q;q^2)_{n}(-q^2;q^2)_{n}q^n}{(-q^2;q^2)_{n}(-q^3;q^2)_{n}} \nonumber \\
&= \frac{1}{1+q}\sum_{n=0}^{\infty}\frac{(q;q^2)_{n}q^n}{(-q^3;q^2)_{n}} \nonumber \\
&= \frac{1}{1+q}{_2}\phi_{1}\bigg(\genfrac{}{}{0pt}{} {q^2,  q}{ -q^3}; q^2,q\bigg).
\end{align}
Replacing $q$ by $q^2$ and setting $(a,b,c,z)=(q^2,q,-q^3,q)$ in \eqref{gasperp15}, we deduce that
\begin{align}\label{4-continue}
{_2}\phi_{1}\bigg(\genfrac{}{}{0pt}{} {q^2,  q} {-q^3}; q^2,q\bigg)&= \frac{(q^3;q^2)_{\infty}}{(q;q^2)_{\infty}}\sum_{n=0}^{\infty}\frac{(q^2,-q^2;q^2)_{n}}{(q^2,-q^3,q^3;q^2)_{n}}\cdot (-q^2)^{n}q^{n(n-1)} \nonumber\\
&=\frac{1}{1-q} \sum_{n=0}^{\infty}\frac{(-q^2;q^2)_{n}(-1)^{n}q^{n(n+1)}}{(-q^3;q^2)_{n}(q^3;q^2)_{n}} \nonumber \\
&=\frac{1}{1-q} \sum_{n=0}^{\infty}\frac{(-q^2;q^2)_{n}}{(-q^3;q^2)_{n}(q^3;q^2)_{n}}\lim_{a \to 0}(\frac{q^2}{a};q^2)_{n}a^n \nonumber \\
&=\frac{1}{1-q} \lim_{a\to 0}{_3}\phi_{2}\bigg(\genfrac{}{}{0pt}{} {q^2/a, q^2, -q^2}{-q^3, q^3};q^2,a  \bigg).
\end{align}
Replacing $q$ by $q^2$ and setting $(b,c,d,\alpha,\beta)=(1,-q^3,q^3,q^2,-q^2)$ in \eqref{liu2089}, we get
\begin{align}\label{4-middle}
&\frac{(q^4,a;q^2)_{\infty}}{(q^2a,q^2;q^2)_{\infty}}{_3}\phi_{2}\bigg(\genfrac{}{}{0pt}{} {q^2/a, q^2, -q^2}{-q^3, q^3};q^2,a  \bigg) \nonumber \\
&=\sum_{n=0}^{\infty}\frac{(1-q^{4n+2})(q^2,q^2/a,q^2;q^2)_{n}(-a)^{n}q^{n(n-1)}}{(1-q^2)(q^2,q^2a,q^2;q^2)_{n}} \cdot {_3}\phi_{2}\bigg(\genfrac{}{}{0pt}{}{q^{-2n},  q^{2n+2},  -q^2} {-q^3, q^3};q^2,q^2 \bigg).
\end{align}
We can now take the limit $a\to 0$ so that
\begin{align}\label{4-final}
\lim_{a\to 0}{_3}\phi_{2}\bigg(\genfrac{}{}{0pt}{} {q^2/a, q^2, -q^2}{-q^3, q^3};q^2,a  \bigg)=\sum_{n=0}^{\infty}(1-q^{4n+2})q^{2n^2}{_3}\phi_{2}\bigg(\genfrac{}{}{0pt}{}{q^{-2n},  q^{2n+2},  -q^2} {-q^3, q^3};q^2,q^2 \bigg).
\end{align}
Replacing $q$ by $q^2$ and setting $(a,b,c)=(q^2,-q,q)$ in \eqref{gasper221}, we get
\begin{align}\label{4-value}
{_3}\phi_{2}\bigg(\genfrac{}{}{0pt}{} {q^{-2n}, q^{2n+2},  -q^2} {-q^3, q^3};q^2,q^2 \bigg)=\frac{(-q,q;q^2)_{n}}{(-q^3,q^3;q^2)_{n}}(-q^2)^{n}=\frac{1-q^2}{1-q^{4n+2}}(-1)^{n}q^{2n}.
\end{align}
Substituting \eqref{4-value} into \eqref{4-final}, we deduce that
\begin{align}\label{4-end}
\lim_{a\to 0}{_3}\phi_{2}\bigg(\genfrac{}{}{0pt}{} {-q^2,  q^2/a,  q^2} {-q^3, q^3};q^2,a \bigg)=(1-q^2)\sum_{n=0}^{\infty}(-1)^{n}q^{2n^2+2n}.
\end{align}
Substituting \eqref{4-end} into \eqref{4-continue} and then combining with \eqref{4-start}, we complete the proof of \eqref{false-4}.
\end{proof}

\begin{proof}[Proof of \eqref{false-5}]
We observe that
\begin{align}\label{5-start}
\sum_{n=0}^{\infty}\frac{(q;-q)_{n}(-q^2;q^2)_{n}q^n}{(-q;q)_{2n+1}}=\frac{1}{1+q}\sum_{n=0}^{\infty}\frac{(q;-q)_{n}q^{n}}{(-q^3;q^2)_{n}}.
\end{align}
We will next consider the sum on the right with $q$ replaced by $-q$. For this sum we note that
\begin{align}\label{5-continue}
\sum_{n=0}^{\infty}\frac{(-q;q)_{n}(-q)^{n}}{(q^3;q^2)_{n}}=\sum_{n=0}^{\infty}\frac{(-q;q)_{n}(-q)^{n}}{(q^{3/2};q)_{n}(-q^{3/2};q)_{n}}={_3}\phi_{2}\bigg(\genfrac{}{}{0pt}{} {q,  -q, 0}{q^{3/2}, -q^{3/2}}; q,-q\bigg).
\end{align}
Setting $(a,b,c,d,\alpha,\beta)=(1,-1,q^{3/2},-q^{3/2},q^2,0)$ in \eqref{liu2089}, we deduce that
\begin{align}\label{5-middle}
&\frac{(q^3,-q;q)_{\infty}}{(q^2,-q^2;q)_{\infty}}{_3}\phi_{2}\bigg(\genfrac{}{}{0pt}{} {q, -q, 0}{q^{3/2},  -q^{3/2}}; q,-q\bigg) \nonumber \\
&=\sum_{n=0}^{\infty}\frac{1-q^{2n+2}}{1-q^2}\cdot \frac{(q^2,q,-q;q)_{n}q^{n}q^{n(n-1)/2}}{(q,q^2,-q^2;q)_{n}}{_3}\phi_{2}\bigg(\genfrac{}{}{0pt}{} {q^{-n}, q^{n+2},  0}{q^{3/2},  -q^{3/2}};q,q\bigg) \nonumber\\
&=\sum_{n=0}^{\infty}\frac{(1-q^{n+1})q^{n(n+1)/2}}{1-q}{_3}\phi_{2}\bigg(\genfrac{}{}{0pt}{} {q^{-n}, q^{n+2}, 0}{q^{3/2},  -q^{3/2}}; q,q\bigg).
\end{align}
Setting $(b,d,e)=(q^{n+2},q^{3/2},-q^{3/2})$ in \eqref{gasper71-new}, we deduce that
\begin{align}\label{5-gasper}
&{_2}\phi_{1}\bigg(\genfrac{}{}{0pt}{} {q^{-n},  q^{-n-\frac{1}{2}}} {q^{3/2}};q, -q^{n+\frac{3}{2}}\bigg)\nonumber \\
&=(-1)^{n}q^{-\binom{n}{2}}(-q^{3/2};q)_{n}(-q^{3/2})^{-n}{_3}\phi_{2}\bigg(\genfrac{}{}{0pt}{} {q^{-n}, q^{n+2}, 0}{q^{3/2}, -q^{3/2}}; q,q\bigg) \nonumber \\
&= q^{-n(n+2)/2}(-q^{3/2};q)_{n}{_3}\phi_{2}\bigg(\genfrac{}{}{0pt}{} {q^{-n}, q^{n+2}, 0}{q^{3/2}, -q^{3/2}}; q,q \bigg).
\end{align}
Setting $(a,b)=(q^{-n},q^{-n-1/2})$ in \eqref{BD}, we deduce that
\begin{align}\label{5-BD}
{_2}\phi_{1}\bigg(\genfrac{}{}{0pt}{} {q^{-n}, q^{-n-1/2}}{q^{3/2}}; q, -q^{n+\frac{3}{2}}\bigg)
=\frac{(-q;q)_{\infty}(q^{1-n};q^2)_{\infty}(q^{n+3};q^2)_{\infty}}{(q^{3/2};q)_{\infty}(-q^{n+3/2};q)_{\infty}}.
\end{align}
It is clear that if $n\ge 1$ is odd, then the left hand side of \eqref{5-BD} equals 0. Now we assume that $n$ is even. We have
\begin{align}\label{even-value}
&\frac{(-q;q)_{\infty}(q^{1-n};q^2)_{\infty}(q^{n+3};q^2)_{\infty}}{(q^{3/2};q)_{\infty}(-q^{n+3/2};q)_{\infty}} \nonumber \\
&=\frac{(1-q)(-1)^{n/2}q^{-n^2/4}(-q^{3/2};q)_{n}}{1-q^{n+1}}.
\end{align}
Substituting \eqref{5-BD} and \eqref{even-value} into \eqref{5-gasper}, we obtain
\begin{align}\label{5-final}
{_3}\phi_{2}\bigg(\genfrac{}{}{0pt}{} {q^{-n},  q^{n+2}, 0}{q^{3/2}, -q^{3/2}}; q,q\bigg)=\left\{\begin{array}{ll}
0 & \textrm{$n$ is odd}; \\
\frac{1-q}{1-q^{n+1}}(-1)^{n/2}q^{n^2/4+n} & \textrm{$n$ is even}.
\end{array} \right.
\end{align}
Substituting \eqref{5-final} into \eqref{5-middle} and combining with \eqref{5-continue}, we deduce that
\begin{align}\label{5-end}
\frac{1}{1-q}\sum_{n=0}^{\infty}\frac{(-q;q)_{n}(-q)^{n}}{(q^3;q^2)_{n}}&=\sum_{n=0,\, \textrm{$n$ even}}^{\infty}(-1)^{n/2}q^{\frac{3}{4}n^2+\frac{3}{2}n} \nonumber \\
&=\sum_{n=0}^{\infty}(-1)^{n}q^{3n^2+3n}.
\end{align}
Replacing $q$ by $-q$ in \eqref{5-end} and noting that $3n^2+3n=3n(n+1)$ is always even, by \eqref{5-start} we complete our proof.
\end{proof}

\end{document}